\documentclass[12pt,leqno]{article}
\usepackage{amssymb,latexsym,amsmath,amsfonts,graphicx, ebezier}
\usepackage{a4wide}
\usepackage{comment}

\numberwithin{equation}{section}

\newtheorem{theoreme}{Theorem}[section]
\newtheorem{lemme}[theoreme]{Lemma}
\newtheorem{proposition}[theoreme]{Proposition}
\newtheorem{corollaire}[theoreme]{Corollary}

\newcommand{\ds}{\displaystyle}
\newcommand{\DD}{\mathbb D}
\newcommand{\TT}{\mathbb T}
\newcommand{\RR}{\mathbb R}
\newcommand{\NN}{\mathbb N}
\newcommand{\CC}{\mathbb C}
\renewcommand{\Re}{{\rm Re} \,}
\renewcommand{\Im}{{\rm Im} \,}
\newcommand{\cqfd}
 {%
 \mbox{}%
 \nolinebreak%
 \hfill%
 \rule{2mm}{2mm}%
 \medbreak%
 \par%
 }

\begin{document}
\title{Optimal logarithmic estimates in the Hardy-Sobolev space of the disk and stability results}

\author{I. Feki,  H. Nfata and  F. Wielonsky  }
\maketitle
\begin{abstract}
We prove a logarithmic estimate in the Hardy-Sobolev space $H^{k, 2}$, $k$ a positive integer, of the unit disk ${\DD}$. This estimate extends
those previously established by L. Baratchart and M.
Zerner in $H^{1,2}$ 
and by S. Chaabane and I. Feki in $H^{k,\infty}$. We
use it to derive logarithmic stability
results for 
the inverse problem of identifying Robin's coefficients in corrosion detection by
electrostatic boundary measurements and for a recovery interpolation scheme in the Hardy-Sobolev space $H^{k, 2}$ with interpolation points located on the boundary $\TT$ of the unit disk.
\end{abstract}



\section{Introduction}

In this paper, we establish a logarithmic estimate in the Hardy-Sobolev spaces $H^{k, 2}$, $k$ an integer $\geq 1$, of the unit disk $\DD$,
$$H^{k, 2}=\{f\in H^{2},~f^{(j)}\in H^{2},\ j=1,...,k\},$$
where $H^{2}$ denotes the usual Hardy space of
analytic functions in $\DD$ having bounded $L^{2}$
norms on circles of radius $r$ as $r$ tends to 1, and $H^{k,2}$ is
 endowed with the Sobolev norm $\|\cdot\|_{k,2}$ such that
 \begin{equation*} \|f\|_{k,2}^{2}:=\ds\sum_{j=0}^{k}\|f^{(j)}\|_2^{2}\quad\text{with}\quad
 \|f\|_{2} =\left(\frac{1}{2\pi }\int_0^{2\pi}
|f(e^{i\theta})|^{2}d\theta\right)^{1/2},\quad f\in H^{2}.
\end{equation*}
Before stating our result, we review a few results from the literature.
Motivated by system theoretical problems, L. Baratchart and M.
Zerner study in \cite{BZ} an interpolation scheme for analytic functions in $\DD$ from boundary values on the unit circle $\TT$. In particular, they prove an estimate
in the Hardy-Sobolev space $H^{1,2}$ of the disk, which shows that an upper bound on the $L^{2}$ norm on $\TT$ can be derived from the $L^{2}$ norm on a subarc of $\TT$ and that the relation between the two norms is of ${\log \log}/{\log}$-type. In
\cite{AL},  Alessandrini and al. have proved a similar
estimate of $1/\log^{\alpha}$-type, $0<\alpha<1$, in more general planar domains, and with quite different method. They have applied their result to an inverse Robin boundary value problem in corrosion detection by electrostatic boundary measurements.
 For bounded analytic functions in the unit disk $\DD$, S. Chaabane and I. Feki \cite{SI} have established in the uniform norm an estimate of $1/\log$-type in the Hardy-sobolev
spaces $H^{k, \infty}$, for any positive integer $k$, and have shown
that their estimate is optimal.

The case of an annulus $G_{s}={\DD}\backslash {s {\DD}}$, $0<s<1$,
was studied by J. Leblond and al. in \cite{JA}, where they showed a stability result of
${1}/{\log}$-type in
the Hardy-Sobolev space $H^{1, 2}(G_{s})$. Their estimates control  the behavior
of a function with respect to the $L^{2}$ norm on the inner boundary ${s\TT}$ from its
$L^{2}$ norm on the outer boundary ${\TT}$. In the same situation, H.
Meftahi and F. Wielonsky \cite{MW} have given a similar estimate of
${1}/{\log}$-type in the Sobolev spaces $H^{k, 2}(G_{s})$,
$k\geq0$, which, in particular, makes explicit the dependance of this estimate
with respect to the magnitude of the radius $s$ of the
inner boundary. This result was applied, among others, to the geometric inverse problem of estimating the area of an unknown cavity in a bounded planar domain.

In the sequel, the unit ball of
$H^{k, 2}$ will be denoted by
$$B_{k, 2}=\{f\in H^{k, 2}, \|f\|_{k, 2}\leq 1\}.$$
Also, for $I$ a subarc of the unit circle $\TT$ of length $2\pi\lambda$, $\lambda\in(0,1)$,
we set
\begin{equation*}
    \|f\|_{2,I}=\left(\frac{1}{2\pi \lambda}\int_{I}
|f(e^{i\theta})|^{2}d\theta\right)^{1/2},
\end{equation*}
for the $L^{2}$ norm of $f$ on $I$.

We now state our main result.
\begin{theoreme}\label{main-thm}
Let $k\geq 1$ be an integer.
There exists two positive constants
$\alpha_{k}$ and $\gamma_{k}$, depending only on $k$, such that for every $f\in B_{k, 2}$ satisfying
$\|f\|_{2,I}\leq e^{-\gamma_{k}/\lambda}$, $I$ a given subarc of $\TT$ of length $2\pi\lambda$,
we have
\begin{equation}\label{ineq-main}
    \|f\|_{2}\leq\frac{\alpha_{k}}{(\lambda\log(1/\|f\|_{2,I}))^{k}}.
\end{equation}
\end{theoreme}
Note that Theorem \ref{main-thm} is easily extended to any bounded subset of functions in $H^{k,2}$. Clearly, the parameter $\lambda$ in the upper bound can be integrated to the constant $\alpha_{k}$. We prefer to let it in the denominator to show the dependance of the upper bound with respect to the length of the subarc $I$.
Note also that, for the particular case $k=1$, the theorem improves upon \cite[Lemma 4.2]{BZ} since the upper bound in (\ref{ineq-main}) has no log-log term in the numerator. 

Actually, Theorem \ref{main-thm} is optimal as shown by the following proposition.
\begin{proposition}\label{opt}
Assume $I\subset\TT$ is the semi-circle $\{e^{i\theta},-{\pi}/{2}\leq\theta\leq{\pi}/{2}\}$ and for $a>1$, consider the sequence of normalized polynomials in $B_{k,2}$,
$$
f_{n}=u_{n}/\|u_{n}\|_{k,2},\qquad u_{n}(z)=(z-a)^{n},\qquad n>0.
$$
Then the norms $\|f_{n}\|_{2,I}$ on the subarc $I$ tends to 0 as $n$ tends to infinity while 
\begin{equation*}
 0<\beta_{k,a}:=\lim_{n\to\infty}\|f_{n}\|_{2}(\log
(1/\|f_{n}\|_{2,I}))^{k}.
\end{equation*}
Moreover, the limit $\beta_{k,a}$ tends to 1 as the root $a$ of $f_{n}$ tends to infinity.
\end{proposition}
It follows from this proposition that inequality (\ref{ineq-main}) cannot be improved by replacing the constant $\alpha_{k}$ with any function of the norm $\|f\|_{2,I}$ which would tend to zero as $\|f\|_{2,I}$ tends to zero.

In Section \ref{Prelim}, we display preliminary results and in Section \ref{Proofs} we give the proofs of our main results. As an application of our estimates, we
establish in Section \ref{Appli}
a stability result
for the inverse problem of the identification of
Robin's coefficient by boundary measurements. We also use
our result to estimate the rate of convergence of an interpolation scheme for recovering  a function
in $H^{k, 2}$ from its values on a subset of the unit circle ${\TT}$, thus improving
results previously established in \cite{BZ}. For results about the recovery of functions by interpolation schemes in Hardy spaces or in more general settings, the reader can consult \cite{Osi,P1,Tot,Tro} and in particular the monograph \cite{P2}.
\section{Preliminary results}\label{Prelim}
We give in this section some preliminary results which will be
useful for the proof of our main estimates.

For a function $f$ analytic in $\overline{\DD}$, we define the integral means with respect to the $L^{2}$ norm by
$$M_{f}(r)=\left(\ds\frac{1}{2\pi}\int_{0}^{2\pi}|f(re^{i\theta})|^{2}d\theta\right)^{{1}/{2}},\quad 0<r\leq1.$$
Let us recall the well-known Hardy's convexity theorem for the particular case $p=2$, see \cite[p.9]{PD}.
\begin{theoreme}\label{conv}
Let $f$ be analytic in $\overline{\DD}$ and $0<r\leq1$. Then, $\log M_{f}(r)$ is a convex
function of $\log r$ which means that if
\begin{equation*} \log r=\alpha \log r_{1}+(1-\alpha)\log r_{2}\quad\text{with}
    \quad 0<r_{1}<r_{2}\leq1,\quad 0\leq\alpha\leq1,
    \end{equation*}
    then
    \begin{equation*}
    M_{f}(r)\leq [M_{f}(r_{1})]^{\alpha} [M_{f}(r_{2})]^{1-\alpha}.
    \end{equation*}
\end{theoreme}
In the next lemma, we recall an inequality involving the $L^{2}$ means  on $\TT$ of an analytic function and its derivatives of higher-order. It is due to Hardy, Landau and Littlewood, cf. \cite[Theorem 1 and Remark 2.5]{HLL}.
\begin{lemme}Let $1\leq k<n$ be two integers.
There exists a constant $C_{n,k}\geq 1$ such that for all functions $f$ in the Hardy-Sobolev space $H^{n,2}$, having a zero of order $n$ at the origin, we have
\begin{equation}\label{HLL}
    \|f^{(k)}(z)\|_2\leq C_{n,k}\|f(z)\|_2^{1-{k}/{n}}\|f^{(n)}(z)\|_2^{{k}/{n}}.
\end{equation}
\end{lemme}
Note that in inequality (\ref{HLL}) the derivatives are taken with respect to the complex variable $z$ and that the inequality is false without the assumption on the vanishing of $f$ at the origin (consider $f(z)=z^{k}$). There are similar inequalities corresponding to derivatives with respect to the real variable $\theta$, the argument of $z=e^{i\theta}$. These are then the well-known Landau-Kolmogorov type inequalities (in the present case, with respect to the $L^{2}$ norm on the interval $[0,2\pi]$), see \cite{Kup,Kwo,Mit}. These inequalities  are also related to the so-called Gagliardo-Nirenberg interpolation inequalities, see \cite{Gag,Nir}.

Next, we have the following lemma about the mean growth of the derivative of a function analytic in the unit disk, see also \cite[p.80]{PD}.
\begin{lemme}\label{ll} Let $f$ be analytic in $\DD$ and let $0<r<\rho\leq1$. Then
\begin{equation}\label{Minko}
M_{f'}(r)\leq \ds\frac{M_{f}(\rho)}{\rho^{2}-r^{2}}.
\end{equation}
\end{lemme}
\textbf{Proof.}
Let $0<r<\rho<1$ and let $z$ be of modulus $r$. By the Cauchy formula,
$$f'(z)=\frac{1}{2\pi i}\int_{|\zeta|=\rho}\frac{f(\zeta)d\zeta}{(\zeta-z)^{2}}=
\frac{\rho}{2\pi}\int_{0}^{2\pi}
\frac{f(\rho e^{i(t+\theta)})e^{i(t-\theta)}}{(\rho e^{it}-r)^{2}}dt.$$
Making use of the continuous form of the Minkowski's inequality with exponent 2, namely
$$\left[\int\left(\int|h(x,y)|d\nu(x)\right)^{2}d\mu(y)\right]^{1/2}\leq
\int\left[\int|h(x,y)|^{2}d\mu(y)\right]^{1/2}d\nu(x),$$
we deduce that
$$M_{f'}(r)\leq\frac{1}{2\pi}\int_{0}^{2\pi}\frac{M_{f}(\rho)dt}{\rho^{2}-2\rho r\cos t+r^{2}}
=\frac{M_{f}(\rho)}{\rho^{2}-r^{2}}.$$
Since $M_{f}(r)$ increases with $r$, we have
$$M_{f}(1):=\sup_{r\to 1}M_{f}(r)=\lim_{r\to 1}M_{f}(r),$$
and the inequality (\ref{Minko}) is also valid for $\rho=1$.
\cqfd
Referring to the proof of \cite[Lemma 4.1]{BZ} where we note that a factor $1/2$ is missing in the last but one inequality, we get the next
result.
\begin{lemme}\label{l} Let $I$ be a subarc of $\TT$ of length $2\pi\lambda$ and let $f$ be a bounded analytic function in
${\DD}$ such that $\|f\|_{L^{\infty}(\DD)}\leq 1$. Then, for every $z\in
\overline{\DD}$, we have
\begin{equation*} |f(z)|\leq \|f\|_{2,I}^{\lambda(1-|z|)/2}.
\end{equation*}
\end{lemme}
Note that the lemma applies in particular to functions in $H^{k,2}$, $k\geq 1$, since the Hardy-Sobolev spaces are included in the disk algebra of functions analytic in $\DD$ and continuous on $\overline \DD$.

We now prove the following lemma which will be the basis for the proof of our results.
\begin{lemme}\label{M1-r}
Let $k$ be a positive integer and
$f\in H^{k,2}$ such that $M_{f^{(k)}}(1)\leq 1$. Then, for $0<r<1$, we have
\begin{equation}\label{diff-M}
\|f\|_{2}\leq\sum_{s=0}^{k-1}\frac{(1-r)^{s}}{s!}M_{f^{(s)}}(r)+\left[\frac{\log r}{\log M_{f^{(k)}}(r)}\right]^{k}.
\end{equation}
\end{lemme}
\textbf{Proof.} We first consider a function $g$ in $H^{1,2}$ such that $M_{g'}(1)\leq 1$. We write
$$g(z)=\ds\sum_{j=0}^{\infty}a_{j}z^{j}$$
for its series expansion in the unit disk.
From Parseval formula, we
get
\begin{equation}\label{eqbb}
    M_{g}^{2}(r)=\ds\sum_{j=0}^{\infty}|a_{j}|^{2}r^{2j},\quad 0\leq r\leq 1.
\end{equation}
The function $M_{g}^{2}(r)$ is differentiable as a function of $r$.
Differentiating (\ref{eqbb}),
we obtain
\begin{equation*}
2M_{g}'(r)M_{g}(r)=2\ds\sum_{j=0}^{\infty}j|a_{j}|^{2}r^{2j-1}.
\end{equation*}
Now, applying Cauchy-Schwarz inequality on the right-hand side, we
get
\begin{equation}\label{eqcc}
   M_{g}'(r)\leq M_{g'}(r).
\end{equation}
From Theorem \ref{conv} and the assumption that $M_{g'}(1)\leq 1$, we obtain for $0<r<t\leq 1$,
\begin{equation}\label{eqdd}
    M_{g'}(t)\leq M_{g'}(r)^{{\log t}/{\log
    r}}.
\end{equation}
Since
\begin{equation}\label{eqddd}
M_{g}(s)-M_{g}(r)=\ds\int_{r}^{s}M_{g}'(t)dt,
\end{equation}
we derive from
(\ref{eqcc}) and (\ref{eqdd}) that
\begin{equation}\label{Mdiff}
M_{g}(s)-M_{g}(r)\leq
 \frac{\left[t^{{\log M_{g'}(r)}/{\log r}+1}\right]_{r}^{s}}{{\log M_{g'}(r)}/{\log r}+1}\leq
 \frac{\log r}{\log M_{g'}(r)}s^{\log M_{g'}(r)/{\log r}}.
\end{equation}
Now, for $g=f^{(k-1)}\in H^{1,2}$, and $0<r<t\leq 1$, (\ref{Mdiff}) rewrites as
\begin{equation}\label{Mdiff-f}
M_{f^{(k-1)}}(t)\leq M_{f^{(k-1)}}(r)+
 \frac{\log r}{\log M_{f^{(k)}}(r)}t^{\log M_{f^{(k)}}(r)/{\log r}}.
\end{equation}
Writing (\ref{eqddd}) for $f^{(k-2)}$, making use of (\ref{eqcc}), and
integrating both sides of the previous inequality (\ref{Mdiff-f}) with respect to $t$, $0<r\leq t\leq s\leq 1$, we obtain
$$
M_{f^{(k-2)}}(s)-M_{f^{(k-2)}}(r)\leq M_{f^{(k-1)}}(r)(s-r)+
\left(\frac{\log r}{\log M_{f^{(k)}}(r)}\right)^{2}s^{\log M_{f^{(k)}}(r)/{\log r}+1}.
$$
Hence, after one integration, and for $0<r<t\leq 1$, (\ref{Mdiff-f}) leads to
$$
M_{f^{(k-2)}}(t)\leq M_{f^{(k-2)}}(r)+M_{f^{(k-1)}}(r)(t-r)+
\left(\frac{\log r}{\log M_{f^{(k)}}(r)}\right)^{2}t^{\log M_{f^{(k)}}(r)/{\log r}}.
$$
Then, it is easily checked that after performing $(k-2)$ more integrations, and for $t=1$, one ends up with (\ref{diff-M}), which proves the lemma.
\cqfd
\section{Proofs of Theorem \ref{main-thm} and Proposition \ref{opt}}\label{Proofs}
{\bf Proof of Theorem \ref{main-thm}.}
Let $f\in B_{k,2}$, $k\geq 1$, and set
\begin{equation}\label{def-g}
g(z)=z^{k+1}\frac{f(z)}{D_{k}},
\end{equation}
where $D_{k}$ is a constant depending only on $k$, chosen such that
\begin{equation}\label{Dk}
M_{g^{(s)}}(1)\leq 1,\quad s=0,\ldots,k,\quad\text{ and }\quad \|g\|_{\infty}\leq 1.
\end{equation}
Note that such a constant exists. Indeed, on one hand, since we assume $\|f\|_{k,2}\leq 1$ and since
$$(z^{k+1}f(z))^{(s)}=\sum_{j=0}^{s}j!\binom{s}{j}\binom{k+1}{j}z^{k+1-j}f^{(s-j)}(z),
$$
we have
$$\|(z^{k+1}f(z))^{(s)}\|_{2}\leq\sqrt{s+1}\max_{j=0,\ldots,s}j!\binom{s}{j}\binom{k+1}{j}.$$
On the other hand, from the Fejer-Riesz inequality \cite[Theorem 3.13]{PD} applied to $f'$, one can derive that
$$\|f\|_{\infty}\leq\sqrt{\pi+1}\|f\|_{1,2}.$$
Hence, if we take, for instance,
$$D_{k}=\max\left(\sqrt{\pi+1},\sqrt{k+1}\max_{j=0,\ldots,k}j!\binom{k}{j}\binom{k+1}{j}\right),$$
then (\ref{Dk}) is satisfied.

Now, as the function $g$ has a zero of order $k+1$ at the origin,
the Hardy-Landau-Littlewood inequality (\ref{HLL}) on the circle of radius $r$, $r\in(0, 1)$, tells us that
\begin{equation}\label{eqll}
    M_{g^{(k)}}(r)\leq
C_{k+1,k}(M_{g}(r))^{\frac{1}{k+1}}
(M_{g^{(k+1)}}(r))^{\frac{k}{k+1}}.
\end{equation}
Applying Lemma
\ref{ll} to the derivative $g^{(k)}$ with $\rho=1$, we obtain
\begin{equation}\label{eqmm}
M_{g^{(k+1)}}(r)\leq \frac{M_{g^{(k)}}(1)}{1-r^{2}}\leq
\frac{M_{g^{(k)}}(1)}{1-r}.\end{equation}
Moreover, from Lemma \ref{l} and the fact that $\|g\|_{\infty}\leq 1$, we deduce that
\begin{equation}\label{eqjj}
   M_{g}(r)\leq N_{I}^{1-r},\quad N_{I}:=\|g\|_{2,I}^{\lambda/2}\leq 1.
\end{equation}
Hence, plugging (\ref{eqmm}) and (\ref{eqjj}) into (\ref{eqll}) and using
that $M_{g^{(k)}}(1)\leq1$, we obtain
\begin{equation}\label{eqnn}
    M^{k+1}_{g^{(k)}}(r)\leq C_{k+1,k}^{k+1}
    \frac{N_{I}^{1-r}}{(1-r)^{k}}.
    \end{equation}
Let us choose $r$ in such a way that
\begin{equation}\label{choix-rr}
    N_{I}^{1-r}=\frac{1}
    {(\log(1/N_{I}))^{(k+1)C_{k+1,k}}},
\end{equation}
or equivalently
\begin{equation}\label{choix-r}
r=1-(k+1)C_{k+1,k}\frac{\log\log(1/N_{I})}
    {\log(1/N_{I})}.
    \end{equation}
    The right-hand side of the above equation should be less than 1, which is true if $N_{I}$ is smaller than $1/e$. Actually, we want that
\begin{equation}\label{ineq-rr}
    1-\frac2e\leq r<1.
    \end{equation}
    Here, and in the sequel, we assume that
$C_{k+1,k}\geq 2$ which can always do (possibly by weakening inequality (\ref{HLL})). Then, it is easily checked that the first inequality is satisfied if $N_{I}$ is chosen so that
   \begin{equation}\label{cond-N}
   N_{I}\leq e^{-\Gamma_{k}},\quad\text{ where}\quad \Gamma_{k}\geq e
   \quad\text{and}\quad\frac{\log\Gamma_{k}}{\Gamma_{k}}=\frac{2}{e(k+1)C_{k+1,k}},
    \end{equation}
    (note that the real valued function $\log x/x$ is bounded above by $1/e$, maximum value which is attained at $x=e$).
In the sequel, we always assume that $N_{I}$ satisfies (\ref{cond-N}).

From the concavity of the function $\log$, we see that for $1-2/e\leq r<1$, it holds
\begin{equation}\label{ineq-r}
    A(r-1)\leq \log r\leq
r-1,\quad\text{where}\quad
A=-\frac{\log(1-2/e)}{2/e}=1.808....
\end{equation}
Using (\ref{eqnn}), we have
\begin{equation}\label{eqoo}
\frac{\log
M_{g^{(k)}}(r)}{\log r}\geq
\frac{1}{k+1}\frac{\log N_{I}^{({1-r})}}{\log r}+
\frac{\log C_{k+1,k}}{\log r}-
\frac{k}{k+1}\frac{\log (1-r)}{\log
r}.\end{equation}
For the first term in the right hand side of (\ref{eqoo}), we obtain from the first inequality in (\ref{ineq-r}) that
\begin{equation}\label{ineq1}
\frac{\log N_{I}^{({1-r})}}{\log
r}\geq
\frac{(1-r)\log N_{I}}{A(r-1)}=\frac{\log(1/N_{I})}{A}.
\end{equation}
For the last two terms of (\ref{eqoo}), applying the second inequality in
(\ref{ineq-r}), we have
\begin{equation}\label{eqqq}
    \frac{\log C_{k+1,k}}{\log r}-\frac{k}{k+1}\frac{\log (1-r)}{\log
r}\geq
\frac{\log(1-r)}{1-r}-\frac{\log C_{k+1,k}}{1-r},
\end{equation}
and using the value of $r$ in (\ref{choix-r}), we obtain the lower bound
$$
\frac{1}{(k+1)C_{k+1,k}}\left(-\log(1/N_{I})+\frac{\log(1/N_{I})(\log(k+1)+\log\log\log(1/N_{I}))}
{\log\log(1/N_{I})}\right).$$
Since we assume $N_{I}\leq e^{-\Gamma_{k}}$ with $\Gamma_{k}\geq e$, see (\ref{cond-N}), the second term in the parenthesis is positive and we finally get the inequality
\begin{equation}\label{ineq2}
\frac{\log C_{k+1,k}}{\log r}-\frac{k}{k+1}\frac{\log (1-r)}{\log
r}\geq\frac{-\log(1/N_{I})}{(k+1)C_{k+1,k}}.
\end{equation}
Plugging (\ref{ineq1}) and (\ref{ineq2}) into (\ref{eqoo}), we obtain
\begin{equation}\label{ineq3}
\frac{\log M_{g^{(k)}}(r)}{\log r}\geq
\left(\frac1A-\frac{1}{C_{k+1,k}}\right)\frac{\log(1/N_{I})}{k+1},
\end{equation}
where the constant in the parenthesis is positive, recall we have assumed that $C_{k+1,k}\geq 2$ after (\ref{ineq-rr}).

From Lemma \ref{M1-r} and the fact that $M_{g^{(k)}}(1)\leq 1$, we know that
\begin{equation}\label{major-g}
    \|g\|_{2}\leq
\sum_{s=0}^{k-1}\frac{(1-r)^{s}}{s!}M_{g^{(s)}}(r)+
\left[\frac{\log r}{\log M_{g^{(k)}(r)}}
\right]^{k}.
\end{equation}
Inequality (\ref{ineq3}) gives an upper bound for the above bracketed term.
It remains to control the means
$M_{g^{(s)}}(r)$ of the derivatives of orders $s=0,\cdots, k-1$. 
If $k=1$, the sum reduces to the single term $M_{g}(r)$ for which, by (\ref{eqjj}), (\ref{choix-rr}) and the fact that $N_{I}<1/e$, see (\ref{cond-N}), we have
$$M_{g}(r)\leq\frac{1}{(\log(1/N_{I}))^{k}}.$$
We assume now that $k\geq 2$. From the
Hardy-Landau-Littlewood inequality, Theorem \ref{conv} applied with $r_{1}=r^{1/\alpha}$, $r_{2}=1$, $0<\alpha<1$, and Lemma \ref{ll},
    we obtain
\begin{align*}
    M^{s+1}_{g^{(s)}}(r)
    & \leq C^{s+1}_{s+1,s} M_{g}(r)M^{s}_{g^{(s+1)}}(r)
    \leq C^{s+1}_{s+1,s} M_{g}(r)M^{s\alpha}_{g^{(s+1)}}(r_{1})\\
    & \leq C^{s+1}_{s+1,s}M_{g}(r)
\frac{M^{s\alpha}_{g^{(s)}}(r)}{(r^{2}-r_{1}^{2})^{s\alpha}},
\end{align*}
where in the second inequality we have also used (\ref{Dk}). Choosing $\alpha=(k-1)/k$, we get
$$
    M^{1+s/k}_{g^{(s)}}(r)\leq C_{s+1,s}^{s+1}
\frac{M_{g}(r)}{(r^{2}-r_{1}^{2})^{s(k-1)/k}},
$$
implying that
$$M_{g^{(s)}}(r)\leq C_{s+1,s}^{s+1}
\frac{M_{g}(r)^{\frac{1}{1+s/k}}}{(r^{2}-r_{1}^{2})^{s}},$$
where we have used that $C_{s+1,s}\geq 1$ and $r^{2}-r_{1}^{2}<1$. Next, from (\ref{eqjj}), (\ref{choix-rr}), and the inequalities
$$N_{I}<1/e,\qquad1<\frac{C_{k+1,k}}{1+s/k},\quad s=0,\ldots,k-1,$$
we deduce that
$$M_{g}(r)^{\frac{1}{1+s/k}}\leq \frac{1}{(\log(1/N_{I}))^{k}}.$$
Hence,
$$
\sum_{s=0}^{k-1}\frac{(1-r)^{s}}{s!}M_{g^{(s)}}(r)
\leq
\left(\sum_{s=0}^{k-1}\frac{C_{s+1,s}^{s+1}}{s!}\frac{(1-r)^{s}}{(r^{2}-r_{1}^{2})^{s}}\right) \frac{1}{(\log(1/N_{I}))^{k}}.
$$
For the second fraction in the sum, we have
$$\frac{1-r}{r^{2}-r_{1}^{2}}=\frac{1}{r^{2}}\frac{1-r}{1-r^{2/(k-1)}}$$
which is upper bounded by some constant $\widetilde C$ depending only on $k$ since $r$ satisfies the inequalities in (\ref{ineq-rr}). Consequently
\begin{equation}\label{ineq-term1}
\sum_{s=0}^{k-1}\frac{(1-r)^{s}}{s!}M_{g^{(s)}}(r)\leq
\frac{Ce^{C\widetilde C}}
    {(\log(1/N_{I}))^{k}},
\end{equation}
where $C:=\max(C_{s+1,s})$ for $s=0,\ldots,k-1$. 
Making use of (\ref{ineq3}) and (\ref{ineq-term1}) in (\ref{major-g}) we get that there exists an explicit constant $\beta_{k}$ depending only on $k$ such that
$$
\|g\|_{2}\leq\frac{\beta_{k}}{(\log(1/N_{I}))^{k}}.
$$
From the relation (\ref{def-g}) between the functions $f$ and $g$ and the definition of $N_{I}$ in (\ref{eqjj}), we derive that
$$
\|f\|_{2}\leq\frac{2^{k}\beta_{k}D_{k}}{(\lambda\log(D_{k}/\||f\|_{2,I}))^{k}}\leq
\frac{\alpha_{k}}{(\lambda\log(1/\||f\|_{2,I}))^{k}},$$
with $\alpha_{k}=2^{k}\beta_{k}D_{k}$,
which is inequality (\ref{ineq-main}). Finally, observe that the condition (\ref{cond-N}) on $N_{I}$ translates into
$$\|f\|_{2,I}\leq D_{k}e^{-2\Gamma_{k}/\lambda}.$$
Setting $\gamma_{k}=2\Gamma_{k}$, the above inequality is weaker than the assumption made in the theorem.
\cqfd
\noindent{\bf Proof of Proposition \ref{opt}.}
For $a>1$, we consider the sequence of polynomials
\begin{equation*}
u_{n}(z)={(z-a)^{n}}\ ,\quad n\geq 1.
\end{equation*}
We have $\|u_{n}\|_{2}^{2}=I_{n}$ with
$$$$
$$
I_{n}=\frac1{2\pi}\int_{-\pi}^{\pi}(1-2a\cos\theta+a^{2})^{n}d\theta
=(4\pi na)^{-1/2}(1+a)^{2n+1}(1+o(1)),
$$
as $n$ tends to infinity, where in the last equality we have used Laplace method for obtaining the asymptotic estimate, see e.g. \cite[Chapter 3]{Mil}.
For the Sobolev norm of $u_{n}$, we therefore get
\begin{align*}
\|u_{n}\|_{k,2}^{2} &
=I_{n}+n^{2}I_{n-1}+\cdots+n^{2}(n-1)^{2}\dots(n-k+1)^{2}I_{n-k}\\
& = (4\pi na)^{-1/2}n^{2k}(1+a)^{2n-2k+1}(1+o(1)),
\end{align*}
as $n$ tends to infinity.
Let $f_{n}=u_{n}/\|u_{n}\|_{k,2}$ be the normalized function in the Hardy-Sobolev space $H^{k,2}$. Then,
\begin{equation}\label{fnT}
\|f_{n}\|_{2}^{2}=n^{-2k}(1+a)^{2k}(1+o(1)),\qquad\text{as }n\to\infty.
\end{equation}
Moreover,
\begin{equation*}
\|u_{n}\|_{2,I}^{2}=\frac1{\pi}\int_{-\pi/2}^{\pi/2}(1-2a\cos\theta+a^{2})^{n}d\theta
=(\pi na)^{-1}(1+a^{2})^{n+1}(1+o(1)),
\end{equation*}
as $n$ tends to infinity, where in the last equality we have again used Laplace method for obtaining the asymptotic estimate.
This implies that
\begin{equation}\label{fnI}
\|f_{n}\|_{2,I}^{2}=2(\pi na)^{-1/2}n^{-2k}(1+a)^{-2n+2k-1}(1+a^{2})^{n+1}(1+o(1)),\qquad\text{as }n\to\infty,
\end{equation}
which shows in particular that $\|f_{n}\|_{2,I}$ tends to zero as $n$ tends to infinity.
Furthermore, we deduce from (\ref{fnT}) and (\ref{fnI}) that
\begin{equation*}
\|f_{n}\|_{2}\log^{k}(1/\|f_{n}\|_{2,I})=\left(\frac{1+a}{2}\right)^{k}
\log^{k}\left(\frac{(1+a)^{2}}{1+a^{2}}\right)(1+o(1)),\qquad\text{as }n\to\infty,
\end{equation*}
from which the assertions in the proposition follow.
    \cqfd
Let $W^{m,2}(I)$ denote the usual Sobolev space of the subarc $I$, equipped with the norm
$$\|f\|_{W^{m,2}(I)}^{2}=\sum_{k=0}^{m}\|\partial^{k}f/\partial\theta^{k}\|^{2}_{2}.$$
On the unit circle $\TT$, traces of functions of $H^{m,2}$ are linked to functions in $W^{m,2}(\TT)$ by the relation
$$H^{m,2}_{|\TT}=H^{2}_{|\TT}\cap W^{m,2}(\TT),$$
see \cite[Lemma 2]{CFJL}. Let $f\in H^{m,2}$, $m\geq 1$. Then, by \cite[Theorem 3.11]{PD}, $f$ is absolutely continuous and
$$f'(z)=-ie^{-i\theta}\frac{\partial f}{\partial\theta},\quad z=e^{i\theta}.$$
For higher-order derivatives, we then have
$$f^{(j)}(z)=Q^{j}(f)(e^{i\theta}),\quad z=e^{i\theta},\quad 1\leq j\leq m,$$
where $Q$ denotes the differential operator $Q(f)(e^{i\theta})=-ie^{-i\theta}\partial f/\partial\theta$.
Hence, there exists constants $K_{j}$ depending only on $j$ such that
\begin{equation}\label{ineg-HW}
\|f^{(j)}\|_{2,I}\leq K_{j}\|f\|_{W^{j,2}(I)},\quad j\leq m.
\end{equation}
We now state the following corollary of Theorem \ref{main-thm}.
\begin{corollaire}\label{cor}Let $m$ and $k$ be two integers with $0\leq m<k$ and let $\gamma_{k}$ be the constant from Theorem \ref{main-thm}. There exists positive constants $\alpha_{k,m}$ and $\beta_{k,m}$ such
that for every $f\in B_{k, 2}$ with
$\|f\|_{W^{m, 2}(I)}\leq  \beta_{k,m}$,
we have
\begin{equation}\label{ineq-cor}
\|f\|_{m, 2}\leq \frac{\alpha_{k,m}}{(\lambda\log(1/\|f\|_{W^{m, 2}(I)}))^{k-m}}.
\end{equation}
 \end{corollaire}
 {\bf Proof.}
The derivatives $f^{(j)}$ of order $j \in \{0,...,m\}$ belong to the unit ball $B_{k-j, 2}$. For sufficiently small $\beta_{k,m}$, we have, in view of (\ref{ineg-HW}), that the $f^{(j)}$ satisfy the assumptions of Theorem \ref{main-thm}. Hence, there exists positive constants $\alpha_{k- j}$ depending only on $k-j$ such that
\begin{equation}\label{ineq-main-j}
    \|f^{(j)}\|_{2}\leq\frac{\alpha_{k- j}}{(\lambda\log(1/\|f^{(j)}\|_{2,I}))^{k-j}},\quad
    0\leq j\leq m.
\end{equation}
Since
\begin{equation*}
\|f^{(j)}\|_{2,I}\leq K_{j}\|f\|_{W^{m, 2}(I)},
\end{equation*}
we derive from (\ref{ineq-main-j}) and the fact that the map $x \mapsto 1/\log(1/x)$ is increasing with $x$ that there exists constants $\widetilde\alpha_{k- j}$ depending only on $k-j$ such that
\begin{equation*}
\|f^{(j)}\|_{2}\leq\frac{\widetilde\alpha_{k- j}}{(\lambda\log(1/\|f\|_{W^{m, 2}(I)}))^{k-j}}
\leq\frac{\widetilde\alpha_{k- j}}{(\lambda\log(1/\|f\|_{W^{m, 2}(I)}))^{k-m}},
\end{equation*}
where in the second inequality we assume that $\beta_{k,m}\leq e^{-1}$.
By taking squares and summing over all indices $j =0,...,m$, we get a constant $\alpha_{k,m}$ such that (\ref{ineq-cor}) holds true.
\cqfd
\section{Applications to an inverse problem and recovery of functions}\label{Appli}
 In this section, we apply our results to obtain logarithmic stability results
for 
the inverse problem of identifying
Robin's coefficient by boundary measurements.
We also apply our result to find an upper bound on the rate of convergence of a recovery interpolation scheme in $H^{1, 2}$ with points located on a subset of the unit circle ${\TT}$.

For the first inverse problem, we assume that a prescribed flux
$\phi\not\equiv 0$ together with measurements $u_{m}$ are given on a subarc $I$ of $\TT$, and we want to find  a function $q$  on $J={\TT}\backslash I$ such that the solution $u$
of \begin{equation*} (NR)\quad\left\{%
\begin{array}{ccl}
    -\Delta u & =0 & \text{in }\DD,\\
    \partial_{n} u & = \phi & \text{on } I,  \\
    \partial_{n} u+qu & = 0 &\text{on }{\TT}\setminus I,
\end{array}%
\right.
\end{equation*}
also satisfies $u_{|I} =u_{m}$.

Let $K$ be a non-empty connected subset of
$J$ such that $\partial J\cap K =\emptyset$. We suppose that $q$
belongs to the class of admissible Robin coefficients
\begin{equation*} {Q}_{ad}=\{ q\in C^{1}_{0}(\overline{J}),
|q^{(k)}|\leq c',\ k=0,1,\text{ and }q\geq c \text{ on }{K}
\},
\end{equation*} where $c$ and $c'$ are two positive
constants.
For $ q \in {Q}_{ad}$, we denote by $ u_{q} $ the
solution of the  Neumann-Robin problem $ (NR)$.

Let $W_{0}^{1,2}(I)$ denote the closure of $C^{1}_{0}(I)$
in $W^{1,2}(I)$, where $C^{1}_{0}(I)$ stands for the subset of $C^{1}(I)$ consisting of functions $f$ that vanish at the boundary $\partial I$ in $\TT$ together with their derivative $\partial f/\partial\theta$. We refer to \cite{SM,SMJ,CFJL} for  the following results.
\begin{lemme}\label{identlemma2}
{\bf(\cite[Theorem 2]{SMJ},\cite[Theorem 2, Lemma 1]{CFJL})}
 Let  $\phi \geq 0\in W_{0}^{1,2}(I)$,  
 $\phi \not\equiv 0 $ and assume that $q \in {Q}_{ad}$ for some constants $c, c' >0$.
 Then the solution $u_q$ of the inverse problem $(NR)$ belongs to $W^{5/2,2}(\DD)$ and its trace on $\TT$ belongs to $W^{2,2}(\TT)$.

Furthermore, there exist positive constants $ \alpha,  \beta $,
such that for every $q \in {Q}_{ad}$,
we have
\begin{eqnarray}
 u_q \geq \alpha > 0 \text{ in }\overline\DD,\quad  \hbox{and}\quad  ||u_{q}||_{W^{2,2}(\TT)} \leq \beta.\nonumber
\end{eqnarray}
\end{lemme}

The next result answers in particular the identifiability issue for
the inverse problem $(NR)$.
\begin{lemme}\label{identlemma}{\bf(\cite[Theorem 1]{SM})}  
Let $\phi$ be given as in the previous lemma. Then, the mapping
\begin{eqnarray}
\begin{array}{lccl}
    F:&{Q}_{ad}& \longrightarrow& L^2(I)\\
      &                q & \longmapsto & {u_{q}}_{|I}
\end{array}\nonumber
\end{eqnarray}
is well defined, continuous and injective.
\end{lemme}
As an application of Theorem \ref{main-thm}, we state the following stability result.
\begin{theoreme} Let $\phi$
be a positive function  in $W^{1, 2}_{0}(I)$, $q_{1},
q_{2}\in {Q}_{ad}$, and let $u_{1}$, $u_{2}$ be the corresponding
solutions to the problem $(NR)$. Then, there exist positive
constants $\eta$ and $\varepsilon<1$, depending on $I, \phi$ and ${Q}_{ad}$, such that
\begin{equation}\label{Robin-bound}
\|q_{1}-q_{2}\|_{2, J}\leq \ds\frac{\eta}{\log(1/\|u_{1}-u_{2}\|_{2, I})},
\end{equation}
provided that $\|u_{1} - u_{2}\|_{W^{1, 2}(I)} < \varepsilon$. 
\end{theoreme}
Note that this result improves upon \cite[Theorem 3]{CFJL} since the upper bound in (\ref{Robin-bound}) has no log-log term in the numerator.\\
{\bf Proof.} 
The proof follows the one of \cite[Theorem 3]{CFJL}, except that we use our Corollary \ref{cor} instead of the similar but weaker result \cite[Corollary 3]{CFJL}.
\cqfd

As a second application, we give a result about an interpolation scheme for recovering a
function in $H^{k, 2}$, $k\geq 1$, from its values at points located in a subarc $I$ of the
boundary $\TT$ of $\DD$. This scheme has been
studied previously in the Hardy-Sobolev spaces $H^{1,2}$ of a disk \cite{BZ} and 
of an annulus \cite{MW}.

Let $S_{n}=\{x_{1}, \cdots,x_{n}\}$ be a
set of $n$ distinct points on $I$. We will say that
$f_{n}\in H^{k, 2}$ interpolates $f\in H^{k, 2}$ on the interpolation set $S_{n}$
if
\begin{equation}\label{interp}
\forall\ i\in\{1, \cdots,n\},\qquad
f_{n}(x_{i})=f(x_{i}).
\end{equation}
We consider a nested sequence of sets
$S_{1}\subset S_{2}\subset\cdots$, with $S:=\bigcup_{n}S_{n}$ such that $\overline{S}=I$.
Condition (\ref{interp}) does not determine $f_{n}$ uniquely. Among all functions in 
$H^{k,2}$ satisfying (\ref{interp}), we choose the unique function $f_{n}$ of minimal norm. It may be characterized by using the orthogonal decomposition $H^{k,2}=Z_{n}\oplus U_{n}$, where $Z_{n}$ denotes the closed subspace of $H^{k,2}$ consisting of functions vanishing on $S_{n}$ and $U_{n}$ denotes the orthogonal complement of $Z_{n}$. Note that $Z_{n}$ is closed since the Hilbert space $H^{k,2}$ has a reproducing kernel and the evaluation maps on $\TT$ are continuous. Then, $f_{n}$ is obtained by projecting $f$ on $U_{n}$, see \cite[Section 2]{BZ} for details.

Our result is a version of \cite[Theorem 4.3]{BZ} in $H^{k,2}$ with an improved error estimate.
\begin{theoreme} Let $I$ be a subarc of $\TT$ of length $2\pi\lambda$, $\lambda\in (0,1)$, and let $f\in B_{k, 2}$. Set
$h_{n}=\sup_{x\in I}d(x, S_{n})$ where $S_{n}$ is an
interpolation set of $n$ distinct points on $I$ and $d(x, S_{n})$ denotes the radial distance from $x\in I$ to the set $S_{n}$. Let $\alpha_{k}$ and $\gamma_{k}$ be the constants from Theorem \ref{main-thm} and assume that $h_{n}\leq e^{-\gamma_{k}/\lambda}$. Then, 
$$\|f-f_{n}\|_{2}\leq \frac{\alpha_{k}}{(\lambda\log(1/h_{n}))^{k}}.$$
\end{theoreme}
{\bf Proof.} The proof is the same as the one of \cite[Theorem 4.3]{BZ}, except that we use our Theorem \ref{main-thm} instead of \cite[Lemma 4.2]{BZ}.
\cqfd
\noindent{\bf Acknowledgement.} The authors would like to thank the reviewer for his careful reading of the manuscript.
\bibliographystyle{plain}

\obeylines \texttt{
Imed Feki, imed.feki@fss.rnu.tn
Laboratoire LAMHA -  LR 11ES52
Universit\'e de Sfax
Facult\'e des sciences de Sfax
BP 1171, 3018 Sfax, TUNISIE
\medskip
Houda Nfata, houda\_nfata@yahoo.com
Laboratoire LAMHA -  LR 11ES52
Universit\'e de Gafsa
Institut Pr\'eparatoire aux Etudes des Ing\'enieurs
Campus Universitaire Sidi Ahmed Zarrouk  2112 Gafsa, TUNISIE
\medskip
Franck Wielonsky, wielonsky@cmi.univ-mrs.fr
Laboratoire LATP - UMR CNRS 6632
Universit\'e d'Aix-Marseille
CMI 39 Rue Joliot Curie
F-13453 Marseille Cedex 20, FRANCE
}
\end{document}